\documentclass[a4paper,11pt]{amsart}
\usepackage[left=1.5in,right=1.5in,top=1.4in,bottom=1.4in]{geometry}
\usepackage{amsthm,amsmath,amssymb}
\usepackage{cmtt}
\usepackage{enumitem}
\usepackage{graphicx}

\let\doendproof\endproof
\renewcommand\endproof{~\hfill\qed\doendproof}

\long\def\ignore#1{}

\graphicspath{{.},{figs/}}

\newtheorem{theorem}{Theorem}
\newtheorem{lemma}{Lemma}[section]

\newtheorem{corollary}[theorem]{Corollary}

\newtheorem{claim}[lemma]{Claim}

\def\qed{\ifvmode\mbox{ }\else\unskip\fi\hskip 1em plus 10fill$\Box$}

\def\F{{\mathcal F}}

\def\setR{\mathbb{R}}

\let\leq\leqslant
\let\geq\geqslant

\DeclareMathOperator{\col}{col}

\newcommand{\bbox}{\vrule height7pt width4pt depth1pt}

\begin{document}
\title{Note on the number of edges in families with linear union-complexity}

\author[P.~Micek]{Piotr Micek}
\address[P.~Micek]{Theoretical Computer Science Department, Faculty of Mathematics and Computer Science, Jagiellonian University, Krak\'{o}w, Poland}
\email{\mtt piotr.micek@tcs.uj.edu.pl}
\author[R.~Pinchasi]{Rom Pinchasi}
\address[R.~Pinchasi]{Mathematics Department, Technion, Haifa, Israel}
\email{\mtt room@math.technion.ac.il}

\begin{abstract}
We give a simple argument showing that the number of edges in the intersection graph $G$ of a family of $n$ sets in the plane with a linear union-complexity is $O(\omega(G)n)$.
In particular, we prove $\chi(G)\leq \col(G)< 19\omega(G)$ for intersection graph $G$ of a family of pseudo-discs, which improves a previous bound.
\end{abstract}

\maketitle
\section{Introduction}

The maximum size of a clique, i.e., a set of pairwise adjacent vertices, in a graph $G$ is called the \emph{clique number} of $G$ and denoted by $\omega(G)$.
A \emph{proper coloring} of a graph is an assignment of colors to the vertices of the graph such that no two adjacent vertices are assigned the same color.
The minimum number of colors sufficient to color a graph $G$ properly is called the \emph{chromatic number} of $G$ and is denoted by $\chi(G)$.
The \emph{coloring number} of a graph $G$, denoted by $\col(G)$, is a minimum integer $k$ such that there is a linear order of vertices of $G$ such that each vertex has less than $k$ backwards neighbors. 
Some authors prefer the notion of \emph{degeneracy} of $G$ which is simply $\col(G)-1$. 
Clearly, $\omega(G) \leq \chi(G) \leq \col(G)$. 

These three graph parameters can be arbitrarily far apart. 
Complete bipartite graphs have chromatic number $2$ and arbitrarily large coloring number. 
There are various constructions of graphs that are triangle-free (have clique number $2$) and still have arbitrarily large chromatic number. 
The first one was given in 1949 by Zykov \cite{Zyk49}, and the one perhaps best known is due to Mycielski \cite{Myc55}. 

In this paper, we restrict our attention to the relation between these three parameters within intersection graphs of geometric objects in the plane.
The \emph{intersection graph} of a family of sets $\F$ is the graph with vertex set $\F$ and edge set consisting of all pairs of intersecting elements of $\F$.
With a slight abuse of notation we identify the family $\F$ with its intersection graph and use terms such as chromatic number and others referring directly to $\F$.

The study of the chromatic number of families of geometric objects in the plane was initiated in the seminal paper of Asplund and Gr{\"u}nbaum \cite{AG60}, where they proved that every family $\F$ of axis-aligned rectangles satisfies $\chi(\F) \leq 4\omega(\F)^2-3\omega(\F)$. 
This was later improved by Hendler \cite{Hen98} to $\chi(\F)\leq 3\omega(\F)^2-2\omega(\F)-1$.
No construction of families of rectangles with $\chi$ superlinear in terms of $\omega$ is known. 
On the other hand, Burling \cite{Bur65} showed that intersection graphs of axis-aligned boxes in $\setR^3$ with clique number $2$ can have arbitrarily large chromatic number.

More recently, Pawlik \emph{et al.}\ \cite{PKKLMTW-segments} presented a construction of triangle-free families of segments in the plane with arbitrarily large chromatic number. 
Suk \cite{Suk+} proved that families of unit-length segments in the plane have chromatic number bounded by a double exponential function of their clique number. 

Kim, Kostochka and Nakprasit \cite{KKN04} showed that every family $\F$ of translates of a fixed convex compact set satisfies $\col(\F) \leq 3\omega(\F)-2$ and presented families of unit-discs witnessing that this bound is tight.
They also proved that every family $\F$ of homothetic (uniformly scaled) copies of a fixed convex compact set in the plane satisfies $\col(\F) \leq 6\omega(\F)-6$. 
Finally, Aloupis~\emph{et al.}~\cite[Lemma~7]{ACCLS09} showed that for every family $\F$ of pseudo-discs $\col(\F)\leq 48\omega(\F)$. We recall that 
a family $\F$ of homothetic copies of a fixed convex compact set in the plane 
is a special case of a family of pseudo-discs (see below).

In this paper we generalize the latter results. 
Let $\F$ be a family of geometric sets in the plane. 
We assume that $\F$ is nice looking, that is every member in $\F$ is bounded by a simple closed Jordan curve and all these curves are in general position,
i.e.\ any
two of them cross only in finite number of points, where two curves $\alpha,\beta$ are said to cross each other at a point, if $\alpha$ passes from one side of $\beta$ to the other at this point. Moreover, no two curves touch or overlap each other, and no three curves pass through a common point.
The \emph{union complexity} of $\F$ is the number of intersection points of boundaries of two (or more) objects in $\F$ that lie on the boundary of the union of all objects 
in $\F$.
Equivalently, it is the number of boundary segments of members in $\F$ composing the boundary of the union of all members in $\F$.
The main contribution of this paper is the following theorem. 
\begin{theorem}\label{thm:main}
Let $c>0$ be a constant and $\F$ be a nice family of sets in the plane such that any subfamily $\F'\subset\F$ has union complexity at most $c|\F'|$. 
Then the intersection graph of $\F$ has at most 
\[
\big(({ce}/{2}+1)\omega(\F)-1\big)|\F| \ \text{edges.}
\]
\end{theorem}

\begin{figure}[ht]
\includegraphics[scale=0.5]{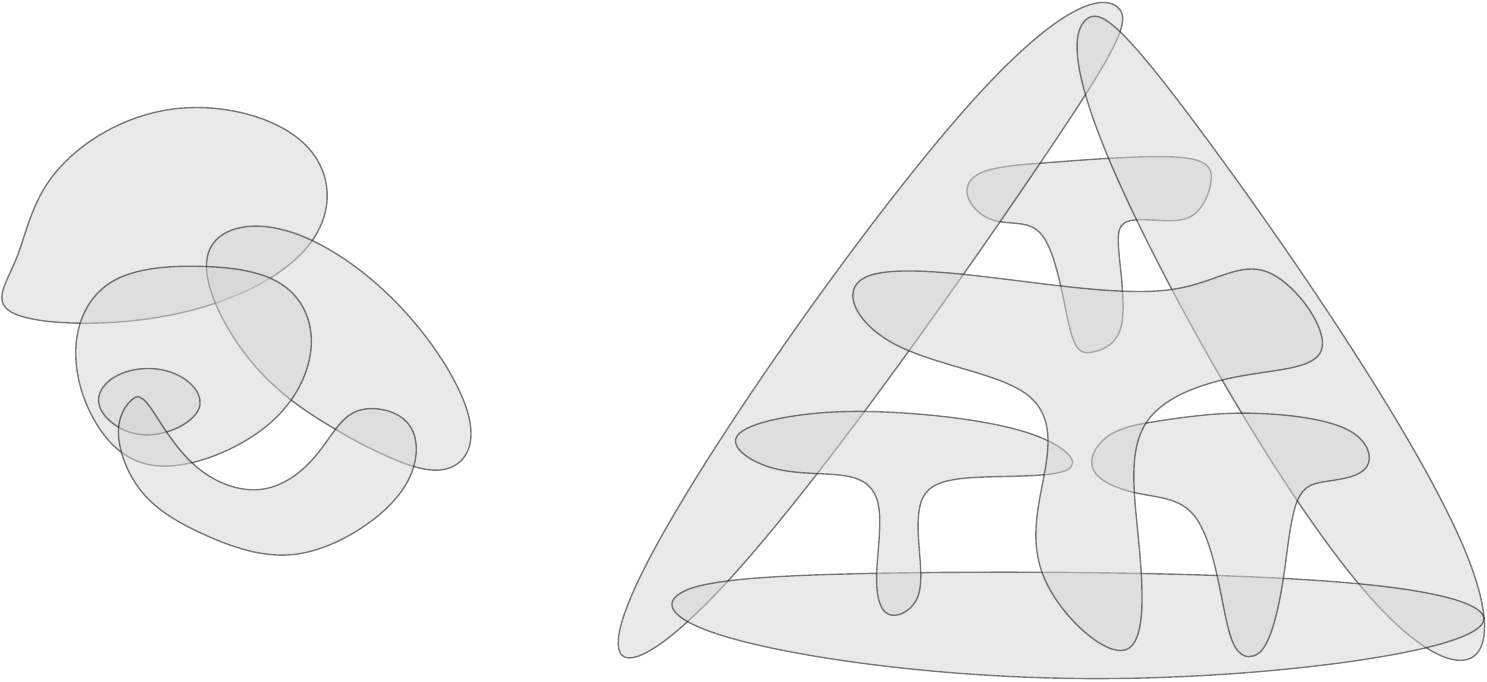}
\caption{On the left a family of pseudo-discs with clique number $3$. On the right, a family of $n=7$ pseudo-discs with union-complexity $6n-12$.}
\label{figure:tight-example}
\end{figure}

Recall that a family $\F$ of nice sets is a family of \emph{pseudo-discs} 
if the boundaries of every two sets from $\F$ intersect in at most two points (i.e.\ in $0$ or $2$ points). 
It is well known, as well as a good exercise to show, 
that a nice family of homothethic copies of a fixed compact convex set is a family of pseudo-discs. 
Kedem \emph{et al.}\ \cite{KLPS86} have shown that the union complexity of a family $\F$ of pseudo-discs, with $|\F|\geq3$, is bounded by $6|\F|-12$ and in fact this bound is best 
possible (see Figure \ref{figure:tight-example}). 
For more examples of families with linear union complexity consult the excellent survey paper \cite{APS08}.

\begin{corollary}\label{cor:pseudo-discs}
A family $\F$ of $n$ pseudo-discs has at most $((3e+1)\omega(\F)-1)n$ edges. 
Therefore
\[
\chi(\F) \leq \col(\F) < (6e+2)\omega(\F) < 19 \omega(\F).
\]
\end{corollary}

This improves the bound given in~\cite{ACCLS09}.
We remark that like the authors of~\cite{ACCLS09} we use probabilistic sampling technique which is usually attributed to Clarkson and Shor~\cite{CS89}.

Note that a bound on the number of edges in a family $\F$ of pseudo-discs in terms of $\omega(\F)$ can be derived from the following result of Fox and Pach in~\cite{FP10}: 
An intersection graph of family of arc-connected sets in the plane with no subgraph isomorphic to $K_{t,t}$ has number of edges bounded in terms of $t$. 
However, the bounding function, obtained this way in \cite{FP10} is quadratic in $\omega(\F)$. 

In a parallelly developing manuscript, Pach and Walczak \cite{PW+} consider families $\F$ of nice-looking sets such that any subfamily $\F'\subset\F$ has union complexity $o(|\F'|^2)$. 
They also give a bound on the number of edges in $\F$ of the form  $f(\omega(\F))|\F|$ but, again, their bounding function $f$ is very large.

Finally, note that one cannot hope for a similar bound, as in Corollary 
\ref{cor:pseudo-discs}, on the number
of edges in the intersection graph of a family of \emph{pseudo-circles} in the plane. 
A family of pseudo-circles is a family of simple closed Jordan curves such that any two intersect in at most two points. 
Indeed, one can easily represent any complete bipartite graph as an intersection graph of pseudo-circles. 
Moreover, Pawlik \emph{et al.}\ \cite{PKKLMTW-shapes} constructed triangle-free families of circles (or eg.\ axis-aligned square frames) with arbitrarily large chromatic number.

In Section \ref{sec:proof} we prove Theorem \ref{thm:main}. 
In Section \ref{sec:partial} we present partial result towards improving the constant in Corollary \ref{cor:pseudo-discs}. 

The crucial point in the proof of Theorem \ref{thm:main} is to bound from above the number of intersection points of boundaries of sets (pseudo-discs) from $\F$ that belong to at most $k$ sets (pseudo-discs) from $\F$.
Given a collection of pseudo-discs $\F$, we denote by $g(\F,k)$ the number of intersection points of two boundaries of sets in $\F$ that are contained in at most $k$ sets from $\F$.
Notice that the result in \cite{KLPS86} implies $g(\F,2)\leq 6n-12$ for 
every collection $\F$ of pseudo-discs. This bound is best possible.
It will follow from the proof of Theorem \ref{thm:main} that  
$g(\F,k)\leq 3ekn$. 
This answer is most likely not best possible in terms of the multiplicative constant. 
It is not hard to see that $O(kn)$ is indeed the correct 
order of magnitude for this question. The fight is for the multiplicative 
constant in front of the term $kn$. 
Improving this constant will directly improve 
the multiplicative constant in 
Corollary \ref{cor:pseudo-discs} and in the best case it can possibly
match and even improve the bound for the special case of a family $\F$ of 
homothetic 
copies of a fixed convex compact set in the plane studied in \cite{KKN04}.
Because improving a multiplicative constant that is already not very far from
optimal can be a very delicate issue, we will be interested at this point
only in tight bounds that we can derive (for some restricted cases).
In the next theorem we obtain a tight upper bound on $g(\F,k)$ in the special 
case where $\F$ is a family of pseudo-discs containing a common point.

\begin{theorem}\label{theorem:pseudo1}
Let $\F$ be a family of $n$ pseudo-discs all containing a common point $O$
in their interior. Then for every $k$ we have $g(\F,k) \leq 2(k-1)n$.
\end{theorem}

\section{Proof of Theorem \ref{thm:main}}\label{sec:proof}

Let $|\F|=n$. 
Let $Z$ be the set of intersection points of boundaries of sets in $\F$. 
We first aim to give an upper bound for $|Z|$.

Pick every member in $\F$ with probability $p$, to be determined later. 
Denote by $\F^*$ the family of those sets in $\F$ that were picked and let $n^*=|\F^*|$. 
Let $S^*$ denote the random variable which is the number of intersection points on the boundary of the union of the sets in $\F^*$. 
By our assumption on the union-complexity of subfamilies of $\F$ we have $S^* \leq cn^*$.
Taking expectations we get 
\[
E(S^*) \leq E(cn^*) = cpn.
\]

In order to estimate $E(S^*)$ from below consider an intersection point $x\in Z$ of the boundaries of two sets $A$ and $B$ in $\F$. 
Observe that $x$ belongs to the boundary of the union of the sets in $F^*$ if and only if both $A$ and $B$ were picked to be in $\F^*$ and every other set in $\F$ that contains $x$ was not picked.
Clearly, the number of those sets in $\F$, different from $A$ and from $B$, that contain $x$ is at most $\omega(\F)-2$. 
Therefore, $x$ appears on the boundary of the union of all sets in $\F^*$ with probability of at least $p^2(1-p)^{\omega(\F)-2}$. 
This gives us a lower bound $|Z|p^2(1-p)^{\omega(\F)-2} \leq E(S^*)$.
Hence
\[
|Z|p^2(1-p)^{\omega(\F)-2} \leq cpn.
\]

We take $p=\frac{1}{\omega(\F)}$ to obtain $|Z| \leq ce\omega(\F)n$, where $e$ denotes the basis of the natural logarithm.
Because every two boundaries of sets in $\F$ intersect zero or two times, we conclude that the number of those pairs in $\F$ whose boundaries intersect is at most $\frac{ce}{2}\omega(\F)n$.

Notice that if two sets in $\F$ intersect, then either their boundaries intersect or one of them is contained in the other.
Clearly, every set in $\F$ is contained in at most $\omega(\F)-1$ other sets from $\F$.
Thus, there can be at most $(\omega(\F)-1)n$ intersections of sets from $\F$ with one set being contained in the other. 
Now when two boundaries of sets from $\F$ intersect by the fact that $\F$ is nice we know that they intersect in at least two points.
It follows that the intersection graph of $\F$ consists of at most $((\frac{ce}{2}+1)w(\F)-1)n$ edges.
\bbox

\section{Towards improving the constant in Theorem \ref{thm:main} for pseudo-discs}\label{sec:partial}

\noindent {\bf Proof of Theorem \ref{theorem:pseudo1}.}
By a result of Snoeyink and Hershberger \cite{SH91}, any family of 
pseudo-circles surrounding a common point can be
\emph{swept by a ray}. In other words, it can be realized as a family of
2-intersecting bi-infinite $x$-monotone curves 
(see [9] for the formal definition of a sweeping)
and this can be done by a one to one continuous transformation
of the plane, after we identify
the two ends at infinity of each curve.

Hence, Theorem \ref{theorem:pseudo1} will follow from the following 
equivalent theorem.
We recall that a family of \emph{pseudo-parabolas} is a family of bi-infinite 
$x$-monotone curves every two of which intersect in at most two points.
A collection of graphs of quadratic polynomials is a natural example
for such a family.

\begin{theorem}\label{theorem:pseudo2}
Let $\F$ be a family of $n$ pseudo-parabolas in the plane.
Then the numbers of intersection points of two curves in $\F$
that lie strictly above at most $k-2$ other curves in $\F$ is bounded
from above by $2(k-1)n$. 
\end{theorem}

\noindent {\bf Proof.}
Consider an intersection point $X$ of two curves $p_{1}$ and $p_{2}$ in $\F$ 
that lies above at most 
$k-2$ curves from $\F$. Without loss of generality
assume that in a small neighborhood to the left of $X$ the curve 
$p_{2}$ lies above 
$p_{1}$ (and consequently $p_{1}$ lies above $p_{2}$ in a small neighborhood
to the right of $X$).

If $X$ is the leftmost intersection point of $p_{1}$ and $p_{2}$, then 
we charge $X$ to $p_{2}$ and we say that this charging is colored red.
If $X$ is not the leftmost intersection point of $p_{1}$ and $p_{2}$
(recall that any two curves in $\F$ intersect at most twice), then 
we charge $X$ to $p_{1}$ and we say that this charging is colored blue.

\begin{claim}\label{claim:1}
No curve in $\F$ can be charged more than $k-1$ times in charging that is 
colored red.
\end{claim}

\noindent {\bf Proof.}
Assume to the contrary that a curve $p$ in $\F$ is charged at least $k$
times with a red charging. Let $X_{1}, X_{2}, \ldots, X_{k}$
denote the $k$ leftmost intersection points on $p$ that charged $p$ 
with a red charging, indexed from leftmost to rightmost.
Let $p_{1}, \ldots, p_{k}$ denote those curves in $\F$ that intersect $p$
at $X_{1}, \ldots, X_{k}$, respectively.
By the definition of a red charging, for every $1 \leq j \leq k$
the curve $p$ lies above the curve $p_{j}$ at any place to the left of 
$X_{j}$. This is because $X_{j}$ is the leftmost intersection point of $p_{j}$
and $p$ and in a small neighborhood to the left of $X_{j}$ the curve $p$ lies
above $p_{j}$.
In particular, the point $X_{1}$ lies above all the $k-1$ curves $p_{2}, \ldots, p_{k}$ which is a contradiction. 
\bbox

\bigskip

Similar to Claim \ref{claim:1} we have have the following observation 
whose proof we omit:

\begin{claim}\label{claim:2}
No curve in $\F$ can be charged more than $k-1$ times in charging that is 
colored blue.
\end{claim}


\bigskip

It follows from Claim \ref{claim:1} and Claim \ref{claim:2} that there are
at most $2(k-1)n$ intersection points of curves in $\F$ that lie strictly 
above at most $k-2$ curves of $\F$.
\bbox

\bigskip

\begin{figure}[ht]
\includegraphics[scale=1]{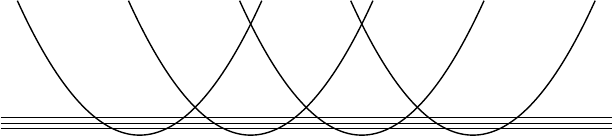}
\caption{Example for $k=4$ and $n=7$. Family contains $2(k-1)(n-k+1)$ intersection points above at most $k-2$ curves.}
\label{figure:partial-example}
\end{figure}

\noindent {\bf Remark.} 
The bound in Theorems \ref{theorem:pseudo1} and \ref{theorem:pseudo2} can indeed be attained
up to a constant additive term that does not depend on $n$.
Consider for instance $k-1$ distinct parallel lines $\ell_{1}, \ldots, \ell_{k-1}$ all of the form $\ell_{i}=\{y=y_i\}$ for $0< y_i < 1/4$. 
Let $p_{1}, \ldots, p_{n-k+1}$ be the parabolas defined by
$p_{i}=\{y=(x-i)^2\}$. Let $\F$ consist of the $n$ curves
$\ell_{1}, \ldots, \ell_{k-1}$ and $p_{1}, \ldots, p_{n-k+1}$. 
Each of the curves $p_{i}$ intersects 
the lines $\ell_{1}, \ldots, \ell_{k-1}$ at $2(k-1)$ intersection points each of
which lies above at most $k-2$ curves of $\F$, that is $\ell_{1}, \ldots, \ell_{k-2}$ (see Figure \ref{figure:partial-example}). 
This gives a count of $2(k-1)(n-k+1)$ such points.

\bibliographystyle{plain}
\bibliography{pseudodiscs}

\begin{thebibliography}{10}

\bibitem{APS08}
Pankaj~K. Agarwal, J{\'a}nos Pach, and Micha Sharir.
\newblock State of the union (of geometric objects).
\newblock In {\em Surveys on discrete and computational geometry}, volume 453
  of {\em Contemp. Math.}, pages 9--48. Amer. Math. Soc., Providence, RI, 2008.

\bibitem{ACCLS09}
Greg Aloupis, Jean Cardinal, S{\'e}bastien Collette, Stefan Langerman, and
  Shakhar Smorodinsky.
\newblock Coloring geometric range spaces.
\newblock {\em Discrete Comput. Geom.}, 41(2):348--362, 2009.

\bibitem{AG60}
Edgar Asplund and Branko Gr\"{u}nbaum.
\newblock On a colouring problem.
\newblock {\em Math. Scand.}, 8:181--188, 1960.

\bibitem{Bur65}
James~P. Burling.
\newblock {\em On coloring problems of families of prototypes}.
\newblock PhD thesis, University of Colorado, Boulder, 1965.

\bibitem{CS89}
Kenneth~L. Clarkson and Peter~W. Shor.
\newblock Applications of random sampling in computational geometry.
\newblock {\em Discrete Comput. Geom.}, 4(5):387--421, 1989.

\bibitem{FP10}
Jacob Fox and J\'{a}nos Pach.
\newblock A separator theorem for string graphs and its applications.
\newblock {\em Combin. Prob. Comput.}, 19(3):371--390, 2010.

\bibitem{Hen98}
Clemens Hendler.
\newblock Schranken f\"{u}r {F}\"{a}rbungs und {C}liquen\"{u}berdeckungszahl
  geometrisch repr\"{a}sentierbarer {G}raphen.
\newblock Master's thesis, Freie Universit\"{a}t Berlin, 1998.

\bibitem{KLPS86}
Klara Kedem, Ron Livn{\'e}, J{\'a}nos Pach, and Micha Sharir.
\newblock On the union of {J}ordan regions and collision-free translational
  motion amidst polygonal obstacles.
\newblock {\em Discrete Comput. Geom.}, 1(1):59--71, 1986.

\bibitem{KKN04}
Seog-Jin Kim, Alexandr Kostochka, and Kittikorn Nakprasit.
\newblock On the chromatic number of intersection graphs of convex sets in the
  plane.
\newblock {\em Electron. J. Combin.}, 11(1):\#R52, 12 pp., 2004.

\bibitem{Myc55}
Jan Mycielski.
\newblock Sur le coloriage des graphes.
\newblock {\em Colloq. Math.}, 3:161--162, 1955.

\bibitem{PW+}
J\'{a}nos Pach and Bartosz Walczak.
\newblock Decomposition of multiple packings with subquadratic union
  complexity.
\newblock submitted.

\bibitem{PKKLMTW-segments}
Arkadiusz Pawlik, Jakub Kozik, Tomasz Krawczyk, Micha\l{} Laso\'{n}, Piotr
  Micek, William~T. Trotter, and Bartosz Walczak.
\newblock Triangle-free intersection graphs of line segments with large
  chromatic number.
\newblock submitted, arXiv:1209.1595.

\bibitem{PKKLMTW-shapes}
Arkadiusz Pawlik, Jakub Kozik, Tomasz Krawczyk, Micha\l{} Laso\'{n}, Piotr
  Micek, William~T. Trotter, and Bartosz Walczak.
\newblock Triangle-free geometric intersection graphs with large chromatic
  number.
\newblock {\em Discrete Comput. Geom.}, 50(3):714--726, 2013.

\bibitem{SH91}
Jack Snoeyink and John Hershberger.
\newblock Sweeping arrangements of curves.
\newblock In {\em Discrete and computational geometry ({N}ew {B}runswick, {NJ},
  1989/1990)}, volume~6 of {\em DIMACS Ser. Discrete Math. Theoret. Comput.
  Sci.}, pages 309--349. Amer. Math. Soc., Providence, RI, 1991.

\bibitem{Suk+}
Andrew Suk.
\newblock Coloring intersection graphs of {$x$}-monotone curves in the plane.
\newblock {\em Combinatorica}.
\newblock to appear, arXiv:1201.0887.

\bibitem{Zyk49}
Alexander~A. Zykov.
\newblock On some properties of linear complexes.
\newblock {\em Mat. Sb. (N.S.)}, 24(66)(2):163--188, 1949.
\newblock in Russian.

\end{thebibliography}

%
%
%
%
%

\end{document}